\documentclass{amsart}
\usepackage{amsfonts}

\newtheorem{theorem}{Theorem}
\theoremstyle{plain}

\newtheorem{corollary}{Corollary}

\newtheorem{lemma}{Lemma}

\newtheorem{proposition}{Proposition}
\newtheorem{remark}{Remark}

\numberwithin{equation}{section}
\input{tcilatex}

\begin{document}
\title[Boas-Bellman Type Inequalities]{Some Boas-Bellman Type Inequalities
in 2-Inner Product Spaces}
\author{S.S. Dragomir}
\address{School of Computer Science and Mathematics\\
Victoria University of Technology\\
PO Box 14428, MCMC \\
Victoria 8001, Australia.}
\email{sever.dragomir@vu.edu.au}
\urladdr{http://rgmia.vu.edu.au/SSDragomirWeb.html}
\author{Y.J. Cho}
\address{Department of Mathematics \\
College of Education\\
Gyeongsang National University\\
Chinju 660-701, Korea}
\email{yjcho@nongae.gsnu.ac.kr}
\author{S.S. Kim$^{\dag }$}
\address{Department of Mathematics\\
Dongeui University\\
Pusan 614-714, Korea.}
\email{sskim@dongeui.ac.kr}
\author{A. Sofo}
\address{School of Computer Science and Mathematics\\
Victoria University of Technology\\
PO Box 14428, MCMC \\
Victoria 8001, Australia.}
\email{sofo@matilda.vu.edu.au}
\date{August 01, 2003.}
\keywords{Bessel's inequality in 2-Inner Product Spaces, Boas-Bellman type
inequalities, 2-Inner Products, 2-Norms.\\
\dag\  \ \ Corresponding author}

\begin{abstract}
Some inequalities in 2-inner product spaces generalizing Bessel's result
that are similar to the Boas-Bellman inequality from inner product spaces,
are given. Applications for determinantal integral inequalities are also
provided.
\end{abstract}

\maketitle

\section{Introduction}

Let $\left( H;\left( \cdot ,\cdot \right) \right) $ be an inner product
space over the real or complex number field $\mathbb{K}$. If $\left(
e_{i}\right) _{1\leq i\leq n}$ are orthonormal vectors in the inner product
space $H,$ i.e., $\left( e_{i},e_{j}\right) =\delta _{ij}$ for all $i,j\in
\left\{ 1,\dots ,n\right\} $ where $\delta _{ij}$ is the Kronecker delta,
then the following inequality is well known in the literature as Bessel's
inequality (see for example \cite[p. 391]{6ba}):%
\begin{equation*}
\sum_{i=1}^{n}\left\vert \left( x,e_{i}\right) \right\vert ^{2}\leq
\left\Vert x\right\Vert ^{2},
\end{equation*}%
for any \ $x\in H.$

For other results related to Bessel's inequality, see \cite{3b} -- \cite{5b}
and Chapter XV in the book \cite{6ba}.

In 1941, R.P. Boas \cite{2b} and in 1944, independently, R. Bellman \cite{1b}
proved the following generalization of Bessel's inequality (see also \cite[%
p. 392]{6ba}).

\begin{theorem}
\label{t1.1}If $x,y_{1},\dots ,y_{n}$ are elements of an inner product space 
$\left( H;\left( \cdot ,\cdot \right) \right) ,$ then the following
inequality: 
\begin{equation*}
\sum_{i=1}^{n}\left\vert \left( x,y_{i}\right) \right\vert ^{2}\leq
\left\Vert x\right\Vert ^{2}\left[ \max_{1\leq i\leq n}\left\Vert
y_{i}\right\Vert ^{2}+\left( \sum_{1\leq i\neq j\leq n}\left\vert \left(
y_{i},y_{j}\right) \right\vert ^{2}\right) ^{\frac{1}{2}}\right] ,
\end{equation*}%
holds.
\end{theorem}

It is the main aim of the present paper to point out the corresponding
version of Boas-Bellman inequality in 2-inner product spaces. Some natural
generalizations and related results are also pointed out. Applications for
determinantal integral inequalities are provided.

For a comprehensive list of fundamental results on 2-inner product spaces
and linear 2-normed spaces, see the recent books \cite{CLKM} and \cite{FC}
where further references are given.

\section{Bessel's Inequality in 2-Inner Product Spaces}

The concepts of $2$-inner products and $2$-inner product spaces have been
intensively studied by many authors in the last three decades. A systematic
presentation of the recent results related to the theory of $2$-inner
product spaces as well as an extensive list of the related references can be
found in the book \cite{CLKM}. Here we give the basic definitions and the
elementary properties of $2$-inner product spaces.

Let $X$ be a linear space of dimension greater than $1$ over the field $%
\mathbb{K}=\mathbb{R}$ of real numbers or the field $\mathbb{K}=\mathbb{C}$
of complex numbers. Suppose that $(\cdot ,\cdot |\cdot )$ is a $\mathbb{K}$%
-valued function defined on $X\times X\times X$ satisfying the following
conditions:

(2I$_{1}$) $(x,x|z)\geq 0$ and $(x,x|z)=0$ if and only if $x$ and $z$ are
linearly dependent,

(2I$_{2}$) $(x,x|z)=(z,z|x),$

(2$I_{3}$) $(y,x|z)=\overline{(x,y|z)},$

(2$I_{4}$) $(\alpha x,y|z)=\alpha (x,y|z)$ for any scalar $\alpha \in 
\mathbb{K},$

($2I_{5}$) $(x+x^{\prime },y|z)=(x,y|z)+(x^{\prime },y|z).$

$(\cdot ,\cdot |\cdot )$ is called a $2$\textit{-inner product}\/ on $X$ and 
$(X,(\cdot ,\cdot |\cdot ))$ is called a $2$\textit{-inner product space\/}
(or $2$\textit{-pre-Hilbert space}). Some basic properties of $2$-inner
product $(\cdot ,\cdot |\cdot )$ can be immediately obtained as follows \cite%
{6b}:

(1) If $\mathbb{K}=\mathbb{R}$, then (2I$_{3}$) reduces to

\begin{equation*}
(y,x|z)=(x,y|z).
\end{equation*}

(2) From (2I$_{3}$) and (2I$_{4}$), we have

\begin{equation*}
(0,y|z)=0,\quad (x,0|z)=0
\end{equation*}

and also

\begin{equation}
(x,\alpha y|z)=\bar{\alpha}(x,y|z).  \label{a1.1}
\end{equation}

(3) Using ($2I_{2}$)--($2I_{5}$), we have

\begin{equation*}
(z,z|x\pm y)=(x\pm y,x\pm y|z)=(x,x|z)+(y,y|z)\pm 2\text{Re}(x,y|z)
\end{equation*}

and

\begin{equation}
\text{Re}(x,y|z)=\frac{1}{4}[(z,z|x+y)-(z,z|x-y)].  \label{a1.2}
\end{equation}

In the real case, (\ref{a1.2}) reduces to

\begin{equation}
(x,y|z)=\frac{1}{4}[(z,z|x+y)-(z,z|x-y)]  \label{a1.3}
\end{equation}

and, using this formula, it is easy to see that, for any $\alpha \in \mathbb{%
R}$, 
\begin{equation}
(x,y|\alpha z)=\alpha ^{2}(x,y|z).  \label{a1.4}
\end{equation}

In the complex case, using (\ref{a1.1}) and (\ref{a1.2}), we have

\begin{equation*}
\text{Im}(x,y|z)=\text{Re}[-i(x,y|z)]=\frac{1}{4}[(z,z|x+iy)-(z,z|x-iy)],
\end{equation*}%
which, in combination with (\ref{a1.2}), yields

\begin{equation}
(x,y|z)=\frac{1}{4}[(z,z|x+y)-(z,z|x-y)]+\frac{i}{4}[(z,z|x+iy)-(z,z|x-iy)].
\label{a1.5}
\end{equation}

Using the above formula and (\ref{a1.1}), we have, for any $\alpha \in 
\mathbb{C}$,

\begin{equation}
(x,y|\alpha z)=|\alpha |^{2}(x,y|z).  \label{a1.6}
\end{equation}

However, for $\alpha \in \mathbb{R}$, (\ref{a1.6}) reduces to (\ref{a1.4}).

Also, from (\ref{a1.6}) it follows that

\begin{equation*}
(x,y|0)=0.
\end{equation*}

(4) For any three given vectors $x,y,z\in X$, consider the vector $%
u=(y,y|z)x-(x,y|z)y.$ By ($2I_{1}$), we know that $(u,u|z)\geq 0$ with the
equality if and only if $u$ and $z$ are linearly dependent. The inequality $%
(u,u|z)\geq 0$ can be rewritten as

\begin{equation}
(y,y|z)[(x,x|z)(y,y|z)-|(x,y|z)|^{2}]\geq 0.  \label{a1.7}
\end{equation}

For $x=z$, (\ref{a1.7}) becomes

\begin{equation*}
-(y,y|z)|(z,y|z)|^{2}\geq 0,
\end{equation*}%
which implies that

\begin{equation}
(z,y|z)=(y,z|z)=0  \label{a1.8}
\end{equation}%
provided $y$ and $z$ are linearly independent. Obviously, when $y$ and $z$
are linearly dependent, (\ref{a1.8}) holds too. Thus (\ref{a1.8}) is true
for any two vectors $y,z\in X.$ Now, if $y$ and $z$ are linearly
independent, then $(y,y|z)>0$ and, from (\ref{a1.7}), it follows that

\begin{equation}
|(x,y|z)|^{2}\leq (x,x|z)(y,y|z).  \label{a1.9}
\end{equation}%
Using (\ref{a1.8}), it is easy to check that (\ref{a1.9}) is trivially
fulfilled when $y$ and $z$ are linearly dependent. Therefore, the inequality
(\ref{a1.9}) holds for any three vectors $x,y,z\in X$ and is strict unless
the vectors $u=(y,y|z)x-(x,y|z)y$ and $z$ are linearly dependent. In fact,
we have the equality in (\ref{a1.9}) if and only if the three vectors $x,y$
and $z$ are linearly dependent.

In any given $2$-inner product space $(X,(\cdot ,\cdot \,|\,\cdot ))$, we
can define a function $\Vert \cdot |\cdot \Vert $ on $X\times X$ by

\begin{equation}
\Vert x|z\Vert =\sqrt{(x,x|z)},  \label{a1.10}
\end{equation}%
for all $x,z\in X$.

It is easy to see that this function satisfies the following conditions:

(2N$_{1}$) $\Vert x|z\Vert \geq 0$ and $\Vert x|z\Vert =0$ if and only if $x$
and $z$ are linearly dependent,

(2N$_{2}$) $\Vert z|x\Vert =\Vert x|z\Vert ,$

(2N$_{3}$) $\Vert \alpha x|z\Vert =|\alpha |\Vert x|z\Vert $ for any scalar $%
\alpha \in \mathbb{K}$,

(2N$_{4}$) $\Vert x+x^{\prime }|z\Vert \leq \Vert x|z\Vert +\Vert x^{\prime
}|z\Vert .$

Any function $\Vert \cdot |\cdot \Vert $ defined on $X\times X$ and
satisfying the conditions ($2N_{1}$)--($2N_{4}$) is called a $2$\textit{-norm%
}\/ on $X$ and $(X,\Vert \cdot |\cdot \Vert )$ is called a \textit{linear }$%
2 $\textit{-normed space \cite{FC}}. Whenever a $2$-inner product space $%
(X,(\cdot ,\cdot |\cdot ))$ is given, we consider it as a linear $2$-normed
space $(X,\Vert \cdot |\cdot \Vert )$ with the $2$-norm defined by (\ref%
{a1.10}).

Let $\left( X;\left( \cdot ,\cdot |\cdot \right) \right) $ be a 2-inner
product space over the real or complex number field $\mathbb{K}$. If $\left(
e_{i}\right) _{1\leq i\leq n}$ are linearly independent vectors in the
2-inner product space $X,$ and, for a given $z\in X,\left(
e_{i},e_{j}|z\right) =\delta _{ij}$ for all $i,j\in \left\{ 1,\dots
,n\right\} $ where $\delta _{ij}$ is the Kronecker delta (we say that the
family $\left( e_{i}\right) _{1\leq i\leq n}$ is $z-$orthonormal), then the
following inequality is the corresponding Bessel's inequality (see for
example \cite{6b}) for $z-$orthonormal family $\left( e_{i}\right) _{1\leq
i\leq n}$ in the 2-inner product space $\left( X;\left( \cdot ,\cdot |\cdot
\right) \right) $:%
\begin{equation}
\sum_{i=1}^{n}\left\vert \left( x,e_{i}|z\right) \right\vert ^{2}\leq
\left\Vert x|z\right\Vert ^{2}\text{,}  \label{1.1}
\end{equation}%
for any \ $x\in X.$ For more details on this inequality, see the recent
paper \cite{6b} and the references therein.

\section{Some Inequalities for 2-Norms}

We start with the following lemma which is also interesting in itself.

\begin{lemma}
\label{l2.1}Let $z_{1},\dots ,z_{n},z\in X$ and $\mu _{1},\dots ,\mu _{n}\in 
\mathbb{K}$. Then one has the inequality: 
\begin{multline}
\left\Vert \sum_{i=1}^{n}\mu _{i}z_{i}|z\right\Vert ^{2}  \label{2.1} \\
\leq \left\{ 
\begin{array}{l}
\max\limits_{1\leq i\leq n}\left\vert \mu _{i}\right\vert
^{2}\sum\limits_{i=1}^{n}\left\Vert z_{i}|z\right\Vert ^{2}; \\ 
\\ 
\left( \sum\limits_{i=1}^{n}\left\vert \mu _{i}\right\vert ^{2\alpha
}\right) ^{\frac{1}{\alpha }}\left( \sum\limits_{i=1}^{n}\left\Vert
z_{i}|z\right\Vert ^{2\beta }\right) ^{\frac{1}{\beta }},\ \ \ \text{where \ 
}\alpha >1,\frac{1}{\alpha }+\frac{1}{\beta }=1; \\ 
\\ 
\sum\limits_{i=1}^{n}\left\vert \mu _{i}\right\vert ^{2}\max\limits_{1\leq
i\leq n}\left\Vert z_{i}|z\right\Vert ^{2},%
\end{array}%
\right. \\
+\left\{ 
\begin{array}{l}
\max\limits_{1\leq i\neq j\leq n}\left\{ \left\vert \mu _{i}\mu
_{j}\right\vert \right\} \sum\limits_{1\leq i\neq j\leq n}\left\vert \left(
z_{i},z_{j}|z\right) \right\vert ; \\ 
\\ 
\left[ \left( \sum\limits_{i=1}^{n}\left\vert \mu _{i}\right\vert ^{\gamma
}\right) ^{2}-\left( \sum\limits_{i=1}^{n}\left\vert \mu _{i}\right\vert
^{2\gamma }\right) \right] ^{\frac{1}{\gamma }}\left( \sum\limits_{1\leq
i\neq j\leq n}\left\vert \left( z_{i},z_{j}|z\right) \right\vert ^{\delta
}\right) ^{\frac{1}{\delta }}, \\ 
\hfill \ \ \ \text{where \ }\gamma >1,\ \ \frac{1}{\gamma }+\frac{1}{\delta }%
=1; \\ 
\\ 
\left[ \left( \sum\limits_{i=1}^{n}\left\vert \mu _{i}\right\vert \right)
^{2}-\sum\limits_{i=1}^{n}\left\vert \mu _{i}\right\vert ^{2}\right]
\max\limits_{1\leq i\neq j\leq n}\left\vert \left( z_{i},z_{j}|z\right)
\right\vert .%
\end{array}%
\right.
\end{multline}
\end{lemma}

\begin{proof}
We observe that 
\begin{align}
\left\Vert \sum_{i=1}^{n}\mu _{i}z_{i}|z\right\Vert ^{2}& =\left(
\sum_{i=1}^{n}\mu _{i}z_{i},\sum_{j=1}^{n}\mu _{j}z_{j}|z\right)  \label{2.2}
\\
& =\sum_{i=1}^{n}\sum_{j=1}^{n}\mu _{i}\overline{\mu _{j}}\left(
z_{i},z_{j}|z\right) =\left\vert \sum_{i=1}^{n}\sum_{j=1}^{n}\mu _{i}%
\overline{\mu _{j}}\left( z_{i},z_{j}|z\right) \right\vert  \notag \\
& \leq \sum_{i=1}^{n}\sum_{j=1}^{n}\left\vert \mu _{i}\right\vert \left\vert 
\overline{\mu _{j}}\right\vert \left\vert \left( z_{i},z_{j}|z\right)
\right\vert  \notag \\
& =\sum_{i=1}^{n}\left\vert \mu _{i}\right\vert ^{2}\left\Vert
z_{i}|z\right\Vert ^{2}+\sum\limits_{1\leq i\neq j\leq n}\left\vert \mu
_{i}\right\vert \left\vert \mu _{j}\right\vert \left\vert \left(
z_{i},z_{j}|z\right) \right\vert .  \notag
\end{align}%
Using H\"{o}lder's inequality, we may write that 
\begin{eqnarray}
&&\sum_{i=1}^{n}\left\vert \mu _{i}\right\vert ^{2}\left\Vert
z_{i}|z\right\Vert ^{2}  \label{2.3} \\
&\leq &\left\{ 
\begin{array}{l}
\max\limits_{1\leq i\leq n}\left\vert \mu _{i}\right\vert
^{2}\sum\limits_{i=1}^{n}\left\Vert z_{i}|z\right\Vert ^{2}; \\ 
\\ 
\left( \sum\limits_{i=1}^{n}\left\vert \mu _{i}\right\vert ^{2\alpha
}\right) ^{\frac{1}{\alpha }}\left( \sum\limits_{i=1}^{n}\left\Vert
z_{i}|z\right\Vert ^{2\beta }\right) ^{\frac{1}{\beta }},\ \ \ \text{where \ 
}\alpha >1,\frac{1}{\alpha }+\frac{1}{\beta }=1; \\ 
\\ 
\sum\limits_{i=1}^{n}\left\vert \mu _{i}\right\vert ^{2}\max\limits_{1\leq
i\leq n}\left\Vert z_{i}|z\right\Vert ^{2}.%
\end{array}%
\right.  \notag
\end{eqnarray}%
By H\"{o}lder's inequality for double sums, we also have 
\begin{equation}
\sum\limits_{1\leq i\neq j\leq n}\left\vert \mu _{i}\right\vert \left\vert
\mu _{j}\right\vert \left\vert \left( z_{i},z_{j}|z\right) \right\vert
\label{2.4}
\end{equation}%
\begin{eqnarray*}
&\leq &\left\{ 
\begin{array}{l}
\max\limits_{1\leq i\neq j\leq n}\left\vert \mu _{i}\mu _{j}\right\vert
\sum\limits_{1\leq i\neq j\leq n}\left\vert \left( z_{i},z_{j}|z\right)
\right\vert ; \\ 
\\ 
\left( \sum\limits_{1\leq i\neq j\leq n}\left\vert \mu _{i}\right\vert
^{\gamma }\left\vert \mu _{j}\right\vert ^{\gamma }\right) ^{\frac{1}{\gamma 
}}\left( \sum\limits_{1\leq i\neq j\leq n}\left\vert \left(
z_{i},z_{j}|z\right) \right\vert ^{\delta }\right) ^{\frac{1}{\delta }}, \\ 
\hfill \ \ \ \text{where \ }\gamma >1,\ \ \frac{1}{\gamma }+\frac{1}{\delta }%
=1; \\ 
\\ 
\sum\limits_{1\leq i\neq j\leq n}\left\vert \mu _{i}\right\vert \left\vert
\mu _{j}\right\vert \max\limits_{1\leq i\neq j\leq n}\left\vert \left(
z_{i},z_{j}|z\right) \right\vert ,%
\end{array}%
\right. \\
&& \\
&=&\left\{ 
\begin{array}{l}
\max\limits_{1\leq i\neq j\leq n}\left\{ \left\vert \mu _{i}\mu
_{j}\right\vert \right\} \sum\limits_{1\leq i\neq j\leq n}\left\vert \left(
z_{i},z_{j}|z\right) \right\vert ; \\ 
\\ 
\left[ \left( \sum\limits_{i=1}^{n}\left\vert \mu _{i}\right\vert ^{\gamma
}\right) ^{2}-\left( \sum\limits_{i=1}^{n}\left\vert \mu _{i}\right\vert
^{2\gamma }\right) \right] ^{\frac{1}{\gamma }}\left( \sum\limits_{1\leq
i\neq j\leq n}\left\vert \left( z_{i},z_{j}|z\right) \right\vert ^{\delta
}\right) ^{\frac{1}{\delta }}, \\ 
\hfill \ \ \ \text{where \ }\gamma >1,\ \ \frac{1}{\gamma }+\frac{1}{\delta }%
=1; \\ 
\\ 
\left[ \left( \sum\limits_{i=1}^{n}\left\vert \mu _{i}\right\vert \right)
^{2}-\sum\limits_{i=1}^{n}\left\vert \mu _{i}\right\vert ^{2}\right]
\max\limits_{1\leq i\neq j\leq n}\left\vert \left( z_{i},z_{j}|z\right)
\right\vert .%
\end{array}%
\right.
\end{eqnarray*}%
Utilizing (\ref{2.3}) and (\ref{2.4}) in (\ref{2.2}), we may deduce the
desired result (\ref{2.1}).
\end{proof}

\begin{remark}
\label{r2.2}Inequality (\ref{2.1}) contains in fact 9 different inequalities
which may be obtained combining the first 3 ones with the last 3 ones.
\end{remark}

A particular result is embodied in the following inequality.

\begin{corollary}
\label{c2.3}With the assumptions in Lemma \ref{l2.1}, we have 
\begin{equation}
\left\Vert \sum_{i=1}^{n}\mu _{i}z_{i}|z\right\Vert ^{2}  \label{2.5}
\end{equation}
\end{corollary}

\begin{align}
& \leq \sum_{i=1}^{n}\left\vert \mu _{i}\right\vert ^{2}\left\{
\max\limits_{1\leq i\leq n}\left\Vert z_{i}|z\right\Vert ^{2}+\frac{\left[
\left( \sum_{i=1}^{n}\left\vert \mu _{i}\right\vert ^{2}\right)
^{2}-\sum_{i=1}^{n}\left\vert \mu _{i}\right\vert ^{4}\right] ^{\frac{1}{2}}%
}{\sum_{i=1}^{n}\left\vert \mu _{i}\right\vert ^{2}}\left(
\sum\limits_{1\leq i\neq j\leq n}\left\vert \left( z_{i},z_{j}|z\right)
\right\vert ^{2}\right) ^{\frac{1}{2}}\right\}  \notag \\
& \leq \sum_{i=1}^{n}\left\vert \mu _{i}\right\vert ^{2}\left\{
\max\limits_{1\leq i\leq n}\left\Vert z_{i}|z\right\Vert ^{2}+\left(
\sum\limits_{1\leq i\neq j\leq n}\left\vert \left( z_{i},z_{j}|z\right)
\right\vert ^{2}\right) ^{\frac{1}{2}}\right\} .  \notag
\end{align}%
The first inequality follows by taking the third branch in the first curly
bracket with the second branch in the second curly bracket for $\gamma
=\delta =2.$

The second inequality in (\ref{2.5}) follows by the fact that 
\begin{equation*}
\left[ \left( \sum_{i=1}^{n}\left\vert \mu _{i}\right\vert ^{2}\right)
^{2}-\sum_{i=1}^{n}\left\vert \mu _{i}\right\vert ^{4}\right] ^{\frac{1}{2}%
}\leq \sum_{i=1}^{n}\left\vert \mu _{i}\right\vert ^{2}.
\end{equation*}%
Applying the following Cauchy-Bunyakovsky-Schwarz inequality 
\begin{equation}
\left( \sum_{i=1}^{n}a_{i}\right) ^{2}\leq n\sum_{i=1}^{n}a_{i}^{2},\ \ \
a_{i}\in \mathbb{R}_{+},\ \ 1\leq i\leq n,  \label{2.6}
\end{equation}%
we may write that 
\begin{equation}
\left( \sum\limits_{i=1}^{n}\left\vert \mu _{i}\right\vert ^{\gamma }\right)
^{2}-\sum\limits_{i=1}^{n}\left\vert \mu _{i}\right\vert ^{2\gamma }\leq
\left( n-1\right) \sum\limits_{i=1}^{n}\left\vert \mu _{i}\right\vert
^{2\gamma }\ \ \ \ \ \ \left( n\geq 1\right)  \label{2.7}
\end{equation}%
and 
\begin{equation}
\left( \sum\limits_{i=1}^{n}\left\vert \mu _{i}\right\vert \right)
^{2}-\sum\limits_{i=1}^{n}\left\vert \mu _{i}\right\vert ^{2}\leq \left(
n-1\right) \sum\limits_{i=1}^{n}\left\vert \mu _{i}\right\vert ^{2}\ \ \ \ \
\ \left( n\geq 1\right) .  \label{2.8}
\end{equation}%
Also, it is obvious that: 
\begin{equation}
\max\limits_{1\leq i\neq j\leq n}\left\{ \left\vert \mu _{i}\mu
_{j}\right\vert \right\} \leq \max\limits_{1\leq i\leq n}\left\vert \mu
_{i}\right\vert ^{2}.  \label{2.9}
\end{equation}%
Consequently, we may state the following coarser upper bounds for $%
\left\Vert \sum_{i=1}^{n}\mu _{i}z_{i}|z\right\Vert ^{2}$ that may be useful
in applications.

\begin{corollary}
\label{c2.4}With the assumptions in Lemma \ref{l2.1}, we have the
inequalities: 
\begin{multline}
\left\Vert \sum_{i=1}^{n}\mu _{i}z_{i}|z\right\Vert ^{2}  \label{2.10} \\
\leq \left\{ 
\begin{array}{l}
\max\limits_{1\leq i\leq n}\left\vert \mu _{i}\right\vert
^{2}\sum\limits_{i=1}^{n}\left\Vert z_{i}|z\right\Vert ^{2}; \\ 
\\ 
\left( \sum\limits_{i=1}^{n}\left\vert \mu _{i}\right\vert ^{2\alpha
}\right) ^{\frac{1}{\alpha }}\left( \sum\limits_{i=1}^{n}\left\Vert
z_{i}|z\right\Vert ^{2\beta }\right) ^{\frac{1}{\beta }},\ \ \ \text{where \ 
}\alpha >1,\frac{1}{\alpha }+\frac{1}{\beta }=1; \\ 
\\ 
\sum\limits_{i=1}^{n}\left\vert \mu _{i}\right\vert ^{2}\max\limits_{1\leq
i\leq n}\left\Vert z_{i}|z\right\Vert ^{2},%
\end{array}%
\right. \\
+\left\{ 
\begin{array}{l}
\max\limits_{1\leq i\leq n}\left\vert \mu _{i}\right\vert
^{2}\sum\limits_{1\leq i\neq j\leq n}\left\vert \left( z_{i},z_{j}|z\right)
\right\vert ; \\ 
\\ 
\left( n-1\right) ^{\frac{1}{\gamma }}\left( \sum\limits_{i=1}^{n}\left\vert
\mu _{i}\right\vert ^{2\gamma }\right) ^{\frac{1}{\gamma }}\left(
\sum\limits_{1\leq i\neq j\leq n}\left\vert z\left( _{i},z_{j}|z\right)
\right\vert ^{\delta }\right) ^{\frac{1}{\delta }}, \\ 
\hfill \ \ \ \text{where \ }\gamma >1,\ \ \frac{1}{\gamma }+\frac{1}{\delta }%
=1; \\ 
\\ 
\left( n-1\right) \sum\limits_{i=1}^{n}\left\vert \mu _{i}\right\vert
^{2}\max\limits_{1\leq i\neq j\leq n}\left\vert \left( z_{i},z_{j}|z\right)
\right\vert .%
\end{array}%
\right.
\end{multline}
\end{corollary}

The proof is obvious by Lemma \ref{l2.1} on applying the inequalities (\ref%
{2.7}) -- (\ref{2.9}).

\begin{remark}
\label{r2.5}The following inequalities which are incorporated in (\ref{2.10}%
) are of special interest: 
\begin{equation}
\left\Vert \sum_{i=1}^{n}\mu _{i}z_{i}|z\right\Vert ^{2}\leq
\max\limits_{1\leq i\leq n}\left\vert \mu _{i}\right\vert ^{2}\left[
\sum\limits_{i=1}^{n}\left\Vert z_{i}|z\right\Vert ^{2}+\sum\limits_{1\leq
i\neq j\leq n}\left\vert \left( z_{i},z_{j}|z\right) \right\vert \right] ;
\label{2.11}
\end{equation}%
\begin{multline}
\left\Vert \sum_{i=1}^{n}\mu _{i}z_{i}|z\right\Vert ^{2}  \label{2.12} \\
\leq \left( \sum\limits_{i=1}^{n}\left\vert \mu _{i}\right\vert ^{2p}\right)
^{\frac{1}{p}}\left[ \left( \sum\limits_{i=1}^{n}\left\Vert
z_{i}|z\right\Vert ^{2q}\right) ^{\frac{1}{q}}+\left( n-1\right) ^{\frac{1}{p%
}}\left( \sum\limits_{1\leq i\neq j\leq n}\left\vert \left(
z_{i},z_{j}|z\right) \right\vert ^{q}\right) ^{\frac{1}{q}}\right] ,
\end{multline}%
where $p>1,$ $\frac{1}{p}+\frac{1}{q}=1;$ and 
\begin{equation}
\left\Vert \sum_{i=1}^{n}\mu _{i}z_{i}|z\right\Vert ^{2}\leq
\sum\limits_{i=1}^{n}\left\vert \mu _{i}\right\vert ^{2}\left[
\max\limits_{1\leq i\leq n}\left\Vert z_{i}|z\right\Vert ^{2}+\left(
n-1\right) \max\limits_{1\leq i\neq j\leq n}\left\vert \left(
z_{i},z_{j}|z\right) \right\vert \right] .  \label{2.13}
\end{equation}
\end{remark}

\section{Some Inequalities for Fourier Coefficients}

The following results holds

\begin{theorem}
\label{t3.1}Let $x,y_{1},\dots ,y_{n},z$ be vectors of a 2-inner product
space $\left( X;\left( \cdot ,\cdot |\cdot \right) \right) $ and $%
c_{1},\dots ,c_{n}\in \mathbb{K}$ $\left( \mathbb{K}=\mathbb{C},\mathbb{R}%
\right) .$ Then one has the inequalities: 
\begin{multline}
\left\vert \sum_{i=1}^{n}c_{i}\left( x,y_{i}|z\right) \right\vert ^{2}
\label{3.1} \\
\leq \left\Vert x|z\right\Vert ^{2}\times \left\{ 
\begin{array}{l}
\max\limits_{1\leq i\leq n}\left\vert c_{i}\right\vert
^{2}\sum\limits_{i=1}^{n}\left\Vert y_{i}|z\right\Vert ^{2}; \\ 
\\ 
\left( \sum\limits_{i=1}^{n}\left\vert c_{i}\right\vert ^{2\alpha }\right) ^{%
\frac{1}{\alpha }}\left( \sum\limits_{i=1}^{n}\left\Vert y_{i}|z\right\Vert
^{2\beta }\right) ^{\frac{1}{\beta }},\ \ \ \text{where \ }\alpha >1,\frac{1%
}{\alpha }+\frac{1}{\beta }=1; \\ 
\\ 
\sum\limits_{i=1}^{n}\left\vert c_{i}\right\vert ^{2}\max\limits_{1\leq
i\leq n}\left\Vert y_{i}|z\right\Vert ^{2};%
\end{array}%
\right.  \\
+\left\Vert x|z\right\Vert ^{2}\times \left\{ 
\begin{array}{l}
\max\limits_{1\leq i\neq j\leq n}\left\{ \left\vert c_{i}c_{j}\right\vert
\right\} \sum\limits_{1\leq i\neq j\leq n}\left\vert \left(
y_{i},y_{j}|z\right) \right\vert ; \\ 
\\ 
\left[ \left( \sum\limits_{i=1}^{n}\left\vert c_{i}\right\vert ^{\gamma
}\right) ^{2}-\left( \sum\limits_{i=1}^{n}\left\vert c_{i}\right\vert
^{2\gamma }\right) \right] ^{\frac{1}{\gamma }}\left( \sum\limits_{1\leq
i\neq j\leq n}\left\vert \left( y_{i},y_{j}|z\right) \right\vert ^{\delta
}\right) ^{\frac{1}{\delta }}, \\ 
\hfill \ \ \ \text{where \ }\gamma >1,\ \ \frac{1}{\gamma }+\frac{1}{\delta }%
=1; \\ 
\\ 
\left[ \left( \sum\limits_{i=1}^{n}\left\vert c_{i}\right\vert \right)
^{2}-\sum\limits_{i=1}^{n}\left\vert c_{i}\right\vert ^{2}\right]
\max\limits_{1\leq i\neq j\leq n}\left\vert \left( y_{i},y_{j}|z\right)
\right\vert .%
\end{array}%
\right. 
\end{multline}
\end{theorem}

\begin{proof}
We note that 
\begin{equation*}
\sum_{i=1}^{n}c_{i}\left( x,y_{i}|z\right) =\left( x,\sum_{i=1}^{n}\overline{%
c_{i}}y_{i}|z\right) .
\end{equation*}%
Using Schwarz's inequality in 2-inner product spaces, we have 
\begin{equation*}
\left\vert \sum_{i=1}^{n}c_{i}\left( x,y_{i}|z\right) \right\vert ^{2}\leq
\left\Vert x|z\right\Vert ^{2}\left\Vert \sum_{i=1}^{n}\overline{c_{i}}%
y_{i}|z\right\Vert ^{2}.
\end{equation*}%
Now using Lemma \ref{l2.1} with $\mu _{i}=\overline{c_{i}},$ $z_{i}=y_{i}$ $%
\left( i=1,\dots ,n\right) ,$ we deduce the desired inequality (\ref{3.1}).
\end{proof}

The following particular inequalities that may be obtained by the
Corollaries \ref{c2.3}, \ref{c2.4}, and Remark \ref{r2.5}, hold.

\begin{corollary}
\label{c3.2}With the assumptions in Theorem \ref{t3.1}, one has the
inequalities: 
\begin{multline}
\left\vert \sum_{i=1}^{n}c_{i}\left( x,y_{i}|z\right) \right\vert ^{2}
\label{3.2} \\
\leq \left\Vert x|z\right\Vert ^{2}\times \left\{ 
\begin{array}{l}
\sum\limits_{i=1}^{n}\left\vert c_{i}\right\vert ^{2}\left\{
\max\limits_{1\leq i\leq n}\left\Vert y_{i}|z\right\Vert ^{2}+\left(
\sum\limits_{1\leq i\neq j\leq n}\left\vert \left( y_{i},y_{j}|z\right)
\right\vert ^{2}\right) ^{\frac{1}{2}}\right\} ; \\ 
\\ 
\max\limits_{1\leq i\leq n}\left\vert c_{i}\right\vert ^{2}\left\{
\sum\limits_{i=1}^{n}\left\Vert y_{i}|z\right\Vert ^{2}+\sum\limits_{1\leq
i\neq j\leq n}\left\vert \left( y_{i},y_{j}|z\right) \right\vert \right\} ;
\\ 
\\ 
\left( \sum\limits_{i=1}^{n}\left\vert c_{i}\right\vert ^{2p}\right) ^{\frac{%
1}{p}}\left\{ \left( \sum\limits_{i=1}^{n}\left\Vert y_{i}|z\right\Vert
^{2q}\right) ^{\frac{1}{q}}+\left( n-1\right) ^{\frac{1}{p}}\left(
\sum\limits_{1\leq i\neq j\leq n}\left\vert \left( y_{i},y_{j}|z\right)
\right\vert ^{q}\right) ^{\frac{1}{q}}\right\} , \\ 
\ \ \text{where \ }p>1,\frac{1}{p}+\frac{1}{q}=1; \\ 
\sum\limits_{i=1}^{n}\left\vert c_{i}\right\vert ^{2}\left\{
\max\limits_{1\leq i\leq n}\left\Vert y_{i}|z\right\Vert ^{2}+\left(
n-1\right) \max\limits_{1\leq i\neq j\leq n}\left\vert \left(
y_{i},y_{j}|z\right) \right\vert \right\} .%
\end{array}%
\right. 
\end{multline}
\end{corollary}

\section{Some Boas-Bellman Type Inequalities in 2-Inner Product Spaces}

If one chooses $c_{i}=\overline{\left( x,y_{i}|z\right) }$ $\left( i=1,\dots
,n\right) $ in (\ref{3.1}), then it is possible to obtain 9 different
inequalities between the Fourier coefficients $\left( x,y_{i}|z\right) $ and
the 2-norms and 2-inner products of the vectors $y_{i}$ $\left( i=1,\dots
,n\right) .$ We restrict ourselves only to those inequalities that may be
obtained from (\ref{3.2}).

From the first inequality in (\ref{3.2}) for $c_{i}=\overline{\left(
x,y_{i}|z\right) },$ we get

\begin{eqnarray*}
&&\left( \sum_{i=1}^{n}\left\vert \left( x,y_{i}|z\right) \right\vert
^{2}\right) ^{2} \\
&\leq &\left\Vert x|z\right\Vert ^{2}\sum_{i=1}^{n}\left\vert \left(
x,y_{i}|z\right) \right\vert ^{2}\left\{ \max\limits_{1\leq i\leq
n}\left\Vert y_{i}|z\right\Vert ^{2}+\left( \sum\limits_{1\leq i\neq j\leq
n}\left\vert \left( y_{i},y_{j}|z\right) \right\vert ^{2}\right) ^{\frac{1}{2%
}}\right\} ,
\end{eqnarray*}%
which is clearly equivalent to the following \textit{Boas-Bellman type
inequality} for 2-inner products:%
\begin{equation}
\sum_{i=1}^{n}\left\vert \left( x,y_{i}|z\right) \right\vert ^{2}\leq
\left\Vert x|z\right\Vert ^{2}\left\{ \max\limits_{1\leq i\leq n}\left\Vert
y_{i}|z\right\Vert ^{2}+\left( \sum\limits_{1\leq i\neq j\leq n}\left\vert
\left( y_{i},y_{j}|z\right) \right\vert ^{2}\right) ^{\frac{1}{2}}\right\} .
\label{4.0}
\end{equation}

From the second inequality in (\ref{3.2}) for $c_{i}=\overline{\left(
x,y_{i}|z\right) ,}$ we get 
\begin{equation*}
\left( \sum_{i=1}^{n}\left\vert \left( x,y_{i}|z\right) \right\vert
^{2}\right) ^{2}\leq \left\Vert x|z\right\Vert ^{2}\max_{1\leq i\leq
n}\left\vert \left( x,y_{i}|z\right) \right\vert ^{2}\left\{
\sum_{i=1}^{n}\left\Vert y_{i}|z\right\Vert ^{2}+\sum_{1\leq i\neq j\leq
n}\left\vert \left( y_{i},y_{j}|z\right) \right\vert \right\} .
\end{equation*}%
Taking the square root in this inequality, we obtain 
\begin{eqnarray}
&&\sum_{i=1}^{n}\left\vert \left( x,y_{i}|z\right) \right\vert ^{2}
\label{4.1} \\
&\leq &\left\Vert x|z\right\Vert \max_{1\leq i\leq n}\left\vert \left(
x,y_{i}|z\right) \right\vert \left\{ \sum_{i=1}^{n}\left\Vert
y_{i}|z\right\Vert ^{2}+\sum_{1\leq i\neq j\leq n}\left\vert \left(
y_{i},y_{j}|z\right) \right\vert \right\} ^{\frac{1}{2}},  \notag
\end{eqnarray}%
for any $x,y_{1},\dots ,y_{n},z$ vectors in the 2-inner product space $%
\left( X;\left( \cdot ,\cdot |\cdot \right) \right) .$

If we assume that $\left( e_{i}\right) _{1\leq i\leq n}$ is an orthonormal
family in $X$ with respect with the vector $z,$ i.e., $\left(
e_{i},e_{j}|z\right) =\delta _{ij}$ for all $i,j\in \left\{ 1,\dots
,n\right\} $, then by (\ref{4.0}) we deduce Bessel's inequality (\ref{1.1}),
while from (\ref{4.1}) we have 
\begin{equation}
\sum_{i=1}^{n}\left\vert \left( x,e_{i}|z\right) \right\vert ^{2}\leq \sqrt{n%
}\left\Vert x|z\right\Vert \max_{1\leq i\leq n}\left\vert \left(
x,e_{i}|z\right) \right\vert ,\ \ \ x\in X.  \label{4.2}
\end{equation}%
From the third inequality in (\ref{3.2}) for $c_{i}=\overline{\left(
x,y_{i}|z\right) ,}$ we deduce 
\begin{multline*}
\left( \sum_{i=1}^{n}\left\vert \left( x,y_{i}|z\right) \right\vert
^{2}\right) ^{2}\leq \left\Vert x|z\right\Vert ^{2}\left(
\sum_{i=1}^{n}\left\vert \left( x,y_{i}|z\right) \right\vert ^{2p}\right) ^{%
\frac{1}{p}} \\
\times \left\{ \left( \sum\limits_{i=1}^{n}\left\Vert y_{i}|z\right\Vert
^{2q}\right) ^{\frac{1}{q}}+\left( n-1\right) ^{\frac{1}{p}}\left(
\sum\limits_{1\leq i\neq j\leq n}\left\vert \left( y_{i},y_{j}|z\right)
\right\vert ^{q}\right) ^{\frac{1}{q}}\right\} ,
\end{multline*}%
for $p>1,$ with $\frac{1}{p}+\frac{1}{q}=1.$ Taking the square root in this
inequality, we get 
\begin{multline}
\sum_{i=1}^{n}\left\vert \left( x,y_{i}|z\right) \right\vert ^{2}\leq
\left\Vert x|z\right\Vert \left( \sum_{i=1}^{n}\left\vert \left(
x,y_{i}|z\right) \right\vert ^{2p}\right) ^{\frac{1}{2p}}  \label{4.3} \\
\times \left\{ \left( \sum\limits_{i=1}^{n}\left\Vert y_{i}|z\right\Vert
^{2q}\right) ^{\frac{1}{q}}+\left( n-1\right) ^{\frac{1}{p}}\left(
\sum\limits_{1\leq i\neq j\leq n}\left\vert \left( y_{i},y_{j}|z\right)
\right\vert ^{q}\right) ^{\frac{1}{q}}\right\} ^{\frac{1}{2}},
\end{multline}%
for any $x,y_{1},\dots ,y_{n},z\in X,$ and $p>1,$ with $\frac{1}{p}+\frac{1}{%
q}=1.$

The above inequality (\ref{4.3}) becomes, for an orthornormal family $\left(
e_{i}\right) _{1\leq i\leq n}$ with respect of the vector $z,$%
\begin{equation}
\sum_{i=1}^{n}\left\vert \left( x,e_{i}|z\right) \right\vert ^{2}\leq n^{%
\frac{1}{q}}\left\Vert x|z\right\Vert \left( \sum_{i=1}^{n}\left\vert \left(
x,e_{i}|z\right) \right\vert ^{2p}\right) ^{\frac{1}{2p}},\ \ \ x\in X.
\label{4.4}
\end{equation}%
Finally, the choice $c_{i}=\overline{\left( x,y_{i}|z\right) }$ $\left(
i=1,\dots ,n\right) $ will produce in the last inequality in (\ref{3.2}) 
\begin{eqnarray*}
&&\left( \sum_{i=1}^{n}\left\vert \left( x,y_{i}|z\right) \right\vert
^{2}\right) ^{2} \\
&\leq &\left\Vert x|z\right\Vert ^{2}\sum_{i=1}^{n}\left\vert \left(
x,y_{i}|z\right) \right\vert ^{2}\left\{ \max\limits_{1\leq i\leq
n}\left\Vert y_{i}|z\right\Vert ^{2}+\left( n-1\right) \max\limits_{1\leq
i\neq j\leq n}\left\vert \left( y_{i},y_{j}|z\right) \right\vert \right\} ,
\end{eqnarray*}%
which gives the following inequality 
\begin{equation}
\sum_{i=1}^{n}\left\vert \left( x,y_{i}|z\right) \right\vert ^{2}\leq
\left\Vert x|z\right\Vert ^{2}\left\{ \max\limits_{1\leq i\leq n}\left\Vert
y_{i}|z\right\Vert ^{2}+\left( n-1\right) \max\limits_{1\leq i\neq j\leq
n}\left\vert \left( y_{i},y_{j}|z\right) \right\vert \right\} ,  \label{4.5}
\end{equation}%
for any $x,y_{1},\dots ,y_{n},z\in X.$

It is obvious that (\ref{4.5}) will give for $z-$orthonormal families, the
Bessel inequality mentioned in $\left( \ref{1.1}\right) $ from Introduction.

\begin{remark}
Observe that, both Boas-Bellman type inequality for 2-inner products
incorporated in (\ref{4.0}) and the inequality (\ref{4.5}) become in the
particular case of $z-$orthonormal families, the regular Bessel's
inequality. Consequently, a comparison of the upper bounds is necessary.

It suffices to consider the quantities%
\begin{equation*}
A_{n}:=\left( \sum\limits_{1\leq i\neq j\leq n}\left\vert \left(
y_{i},y_{j}|z\right) \right\vert ^{2}\right) ^{\frac{1}{2}}
\end{equation*}%
and%
\begin{equation*}
B_{n}:=\left( n-1\right) \max\limits_{1\leq i\neq j\leq n}\left\vert \left(
y_{i},y_{j}|z\right) \right\vert ,
\end{equation*}%
where $n\geq 1,$ and $y_{1},\dots ,y_{n},z\in X.$

If we choose $n=3,$ we have%
\begin{equation*}
A_{3}=\sqrt{2}\left( \left( y_{1},y_{2}|z\right) ^{2}+\left(
y_{2},y_{3}|z\right) ^{2}+\left( y_{3,}y_{1}|z\right) ^{2}\right) ^{1/2}
\end{equation*}%
and%
\begin{equation*}
B_{3}=2\max \left\{ \left\vert \left( y_{1},y_{2}|z\right) \right\vert
,\left\vert \left( y_{2},y_{3}|z\right) \right\vert ,\left\vert \left(
y_{3,}y_{1}|z\right) \right\vert \right\} 
\end{equation*}%
where $y_{1},y_{2},y_{3},z\in X$.

If we consider $a:=\left\vert \left( y_{1},y_{2}|z\right) \right\vert \geq
0,b:=\left\vert \left( y_{2},y_{3}|z\right) \right\vert \geq 0$ and $%
c:=\left\vert \left( y_{3,}y_{1}|z\right) \right\vert \geq 0$, then we have
to compare%
\begin{equation*}
A:=\sqrt{2}\left( a^{2}+b^{2}+c^{2}\right) ^{1/2}
\end{equation*}%
with%
\begin{equation*}
B_{3}=2\max \left\{ a,b,c\right\} .
\end{equation*}%
If we assume that $b=c=1,$ then $A:=\sqrt{2}\left( a^{2}+2\right) ^{1/2},$ $%
B_{3}=2\max \left\{ a,1\right\} .$ Finally, for $a=1,$ we get $A=\sqrt{6}%
,B=2 $ showing that $A>B,$ while for $a=2$ we have $A=\sqrt{12},B=4$ showing
that $B>A.$

In conclusion, we may state that the bounds 
\begin{equation*}
M_{1}:=\left\Vert x|z\right\Vert ^{2}\left\{ \max\limits_{1\leq i\leq
n}\left\Vert y_{i}|z\right\Vert ^{2}+\left( \sum\limits_{1\leq i\neq j\leq
n}\left\vert \left( y_{i},y_{j}|z\right) \right\vert ^{2}\right) ^{\frac{1}{2%
}}\right\}
\end{equation*}%
and%
\begin{equation*}
M_{2}:=\left\Vert x|z\right\Vert ^{2}\left\{ \max\limits_{1\leq i\leq
n}\left\Vert y_{i}|z\right\Vert ^{2}+\left( n-1\right) \max\limits_{1\leq
i\neq j\leq n}\left\vert \left( y_{i},y_{j}|z\right) \right\vert \right\}
\end{equation*}%
for the Bessel's sum $\sum_{i=1}^{n}\left\vert \left( x,y_{i}|z\right)
\right\vert ^{2}$, cannot be compared in general, meaning that some time one
is better than the other.
\end{remark}

\section{Applications for Determinantal Integral Inequalities}

Let $\left( \Omega ,\Sigma ,\mu \right) $ be a measure space consisting of a
set $\Omega ,$ a $\sigma -$algebra $\Sigma $ of subsets of $\Omega $ and a
countably additive and positive measure $\mu $ on $\Sigma $ with values in $%
\mathbb{R\cup \{\infty \}}$.

Denote by $L_{\rho }^{2}\left( \Omega \right) $ the Hilbert space of all
real-valued functions $f$ defined on $\Omega $ that are $2-\rho -$integrable
on $\Omega ,$ i.e., $\int_{\Omega }\rho \left( s\right) \left\vert f\left(
s\right) \right\vert ^{2}d\mu \left( s\right) <\infty ,$ where $\rho :\Omega
\rightarrow \lbrack 0,\infty )$ is a measurable function on $\Omega .$

We can introduce the following 2-inner product on $L_{\rho }^{2}\left(
\Omega \right) $ by formula%
\begin{equation}
\left( f,g|h\right) _{\rho }:=\frac{1}{2}\int_{\Omega }\int_{\Omega }\rho
\left( s\right) \rho \left( t\right) \left\vert 
\begin{array}{cc}
f\left( s\right)  & f\left( t\right)  \\ 
h\left( s\right)  & h\left( t\right) 
\end{array}%
\right\vert \left\vert 
\begin{array}{cc}
g\left( s\right)  & g\left( t\right)  \\ 
h\left( s\right)  & h\left( t\right) 
\end{array}%
\right\vert d\mu \left( s\right) d\mu \left( t\right) ,  \label{5.1}
\end{equation}%
where by%
\begin{equation*}
\left\vert 
\begin{array}{cc}
f\left( s\right)  & f\left( t\right)  \\ 
h\left( s\right)  & h\left( t\right) 
\end{array}%
\right\vert ,
\end{equation*}%
we denote the determinant of the matrix%
\begin{equation*}
\left[ 
\begin{array}{cc}
f\left( s\right)  & f\left( t\right)  \\ 
h\left( s\right)  & h\left( t\right) 
\end{array}%
\right] ,
\end{equation*}%
generating the 2-norm on $L_{\rho }^{2}\left( \Omega \right) $ expressed by%
\begin{equation}
\left\Vert f|h\right\Vert _{\rho }:=\left( \frac{1}{2}\int_{\Omega
}\int_{\Omega }\rho \left( s\right) \rho \left( t\right) \left\vert 
\begin{array}{cc}
f\left( s\right)  & f\left( t\right)  \\ 
h\left( s\right)  & h\left( t\right) 
\end{array}%
\right\vert ^{2}d\mu \left( s\right) d\mu \left( t\right) \right) ^{1/2}.
\label{5.2}
\end{equation}%
A simple calculation with integrals reveals that%
\begin{equation}
\left( f,g|h\right) _{\rho }=\left\vert 
\begin{array}{cc}
\begin{array}{c}
\int_{\Omega }\rho fgd\mu  \\ 
\end{array}
& 
\begin{array}{c}
\int_{\Omega }\rho fhd\mu  \\ 
\end{array}
\\ 
\int_{\Omega }\rho ghd\mu  & \int_{\Omega }\rho h^{2}d\mu 
\end{array}%
\right\vert   \label{5.3}
\end{equation}%
and%
\begin{equation}
\left\Vert f|h\right\Vert _{\rho }=\left\vert 
\begin{array}{cc}
\begin{array}{c}
\int_{\Omega }\rho f^{2}d\mu  \\ 
\end{array}
& 
\begin{array}{c}
\int_{\Omega }\rho fhd\mu  \\ 
\end{array}
\\ 
\int_{\Omega }\rho fhd\mu  & \int_{\Omega }\rho h^{2}d\mu 
\end{array}%
\right\vert ^{1/2},  \label{5.4}
\end{equation}%
where, for simplicity, instead of $\int_{\Omega }\rho \left( s\right)
f\left( s\right) g\left( s\right) d\mu \left( s\right) ,$ we have written $%
\int_{\Omega }\rho fgd\mu .$

Using the representations (\ref{5.3}), (\ref{5.4}) and the inequalities for
2-inner products and 2-norms established in the previous sections, one may
state some interesting determinantal integral inequalities, as follows.

\begin{proposition}
\label{p5.1}Let $f,g_{1},...,g_{n},h\in L_{\rho }^{2}\left( \Omega \right) ,$
where $\rho :\Omega \rightarrow \lbrack 0,\infty )$ is a measurable function
on $\Omega .$ Then we have the inequality%
\begin{eqnarray*}
&&\sum_{i=1}^{n}\left\vert 
\begin{array}{cc}
\begin{array}{c}
\int_{\Omega }\rho fg_{i}d\mu \\ 
\end{array}
& 
\begin{array}{c}
\int_{\Omega }\rho fhd\mu \\ 
\end{array}
\\ 
\int_{\Omega }\rho g_{i}hd\mu & \int_{\Omega }\rho h^{2}d\mu%
\end{array}%
\right\vert ^{2} \\
&\leq &\left\vert 
\begin{array}{cc}
\begin{array}{c}
\int_{\Omega }\rho f^{2}d\mu \\ 
\end{array}
& 
\begin{array}{c}
\int_{\Omega }\rho fhd\mu \\ 
\end{array}
\\ 
\int_{\Omega }\rho fhd\mu & \int_{\Omega }\rho h^{2}d\mu%
\end{array}%
\right\vert \times \left\{ \max_{1\leq i\leq n}\left\vert 
\begin{array}{cc}
\begin{array}{c}
\int_{\Omega }\rho g_{i}^{2}d\mu \\ 
\end{array}
& 
\begin{array}{c}
\int_{\Omega }\rho g_{i}hd\mu \\ 
\end{array}
\\ 
\int_{\Omega }\rho g_{i}hd\mu & \int_{\Omega }\rho h^{2}d\mu%
\end{array}%
\right\vert \right. \\
&& \\
&&\left. +\left( \sum_{1\leq i\neq j\leq n}^{n}\left\vert 
\begin{array}{cc}
\begin{array}{c}
\int_{\Omega }\rho g_{j}g_{i}d\mu \\ 
\end{array}
& 
\begin{array}{c}
\int_{\Omega }\rho g_{j}hd\mu \\ 
\end{array}
\\ 
\int_{\Omega }\rho g_{i}hd\mu & \int_{\Omega }\rho h^{2}d\mu%
\end{array}%
\right\vert ^{2}\right) ^{1/2}\right\} .
\end{eqnarray*}
\end{proposition}

The proof follows by the inequality $\left( \ref{4.0}\right) $ applied for
the 2-inner product and 2-norm defined in $\left( \ref{5.1}\right) $ and $%
\left( \ref{5.2}\right) ,$ and utilizing the identities $\left( \ref{5.3}%
\right) $ and $\left( \ref{5.4}\right) .$

If one uses the inequality (\ref{4.5}), then that one may state the
following result as well

\begin{proposition}
\label{p5.2}Let $f,g_{1},...,g_{n},h\in L_{\rho }^{2}\left( \Omega \right) ,$
where $\rho :\Omega \rightarrow \lbrack 0,\infty )$ is a measurable function
on $\Omega .$ Then we have the inequality%
\begin{eqnarray*}
&&\sum_{i=1}^{n}\left\vert 
\begin{array}{cc}
\begin{array}{c}
\int_{\Omega }\rho fg_{i}d\mu \\ 
\end{array}
& 
\begin{array}{c}
\int_{\Omega }\rho fhd\mu \\ 
\end{array}
\\ 
\int_{\Omega }\rho g_{i}hd\mu & \int_{\Omega }\rho h^{2}d\mu%
\end{array}%
\right\vert ^{2} \\
&\leq &\left\vert 
\begin{array}{cc}
\begin{array}{c}
\int_{\Omega }\rho f^{2}d\mu \\ 
\end{array}
& 
\begin{array}{c}
\int_{\Omega }\rho fhd\mu \\ 
\end{array}
\\ 
\int_{\Omega }\rho fhd\mu & \int_{\Omega }\rho h^{2}d\mu%
\end{array}%
\right\vert \times \left\{ \max_{1\leq i\leq n}\left\vert 
\begin{array}{cc}
\begin{array}{c}
\int_{\Omega }\rho g_{i}^{2}d\mu \\ 
\end{array}
& 
\begin{array}{c}
\int_{\Omega }\rho g_{i}hd\mu \\ 
\end{array}
\\ 
\int_{\Omega }\rho g_{i}hd\mu & \int_{\Omega }\rho h^{2}d\mu%
\end{array}%
\right\vert \right. \\
&& \\
&&\left. +\left( n-1\right) \max_{1\leq i\neq j\leq n}\left\vert 
\begin{array}{cc}
\begin{array}{c}
\int_{\Omega }\rho g_{j}g_{i}d\mu \\ 
\end{array}
& 
\begin{array}{c}
\int_{\Omega }\rho g_{j}hd\mu \\ 
\end{array}
\\ 
\int_{\Omega }\rho g_{i}hd\mu & \int_{\Omega }\rho h^{2}d\mu%
\end{array}%
\right\vert \right\} .
\end{eqnarray*}
\end{proposition}

\textbf{Acknowledgement:}{\small \ S. S. Dragomir and Y. J. Cho greatly
acknowledge the financial support from the Brain Pool Program (2002) of the
Korean Federation of Science and Technology Societies. The research was
performed under the "Memorandum of Understanding" between Victoria
University and Gyeongsang National University.}

\end{document}